\documentclass[10pt]{article}
\usepackage{amsmath}
\usepackage{amsthm}
\usepackage{amssymb}
\usepackage{wasysym}

\newtheorem{proposition}{Proposition}[section]

\newtheorem{lemme}[proposition]{Lemma}
\newtheorem{theorem}{Theorem}[section]

\newdimen\AAdi%
\newbox\AAbo%
\def\AArm{\fam0 }
\def\AAk#1#2{\setbox\AAbo=\hbox{#2}\AAdi=\wd\AAbo\kern#1\AAdi{}}%
\def\AAr#1#2#3{\setbox\AAbo=\hbox{#2}\AAdi=\ht\AAbo\raise#1\AAdi\hbox{#3}
}%
\def\kac{Ka$\mathrm{\check{c}}\ $}        
\def\BBe{{\AArm I\!E}}%
\def\BBone{{\AArm 1\AAk{-.8}{I}I}}%
\def\r{{\mathbf r}}
\newcommand {\CA}{{\cal A}}

\newcommand {\CL}{{\cal L}}

\newcommand {\CO}{{\cal O}}
\newcommand {\CP}{{\cal P}}

\newcommand {\CR}{{\cal R}}

\newcommand{\N}{\mathbb N}
\newcommand{\R}{\mathbb R}
\newcommand{\C}{\mathbb C}
\newcommand{\Z}{\mathbb Z}
\newcommand{\pardef}{\stackrel{\scriptscriptstyle\rm def}{=}}

\newcommand{\8}{\infty}
\newcommand{\disp}{\displaystyle}
\newcommand{\wt}{\widetilde}

\newcommand{\eps}{\varepsilon}

\newcommand{\inte}[1]{\stackrel{\circ}{#1}}
\newcommand{\unli}[1]{\underline{#1}}
\newcommand{\ul}{\unli}

\newcommand{\ninfty}{{n\rightarrow +\8}}

\begin{document}

\title{{\Large Birkhoff averages of Poincar{\'e} cycles}\\ {\Large for Axiom A diffeomorphisms}}
\author{Jean-Ren{\'e} Chazottes$^1$, Renaud Leplaideur$^2$}
\date{\begin{tabular}{c} $^1$ Centre de Physique Th{\'e}orique\\
CNRS-Ecole polytechnique, UMR 7644\\
F-91128 Palaiseau Cedex\\
{\tt jeanrene@cpht.polytechnique.fr}
\end{tabular}
\begin{tabular}{c} $^2$ Laboratoire de Math{\'e}matiques\\
 Universit{\'e} de
Bretagne Occidentale\\ 6, rue Victor Le Gorgeu\\ BP 809 \\
F-29285 Brest Cedex\\
{\tt Renaud.Leplaideur@univ-brest.fr}
\end{tabular}}

\thispagestyle{empty}

\maketitle

\begin{abstract}
We study the time of $n$th return of orbits
to some given (union of) rectangle(s) of a Markov partition 
of an Axiom A diffeomorphism.
Namely, we prove the existence of a
scaled generating function for these returns with respect to
any Gibbs measure (associated to a H{\"o}lderian potential).
As a by-product, we derive precise large deviation estimates and a central
limit theorem for Birkhoff averages of Poincar{\'e}
cycles.
We emphasize that  we look at the limiting behavior in term
of number of visits (the size of the visited set is kept fixed).
Our approach relies on the spectral properties of a one-parameter family of
induced transfer operators on unstable leaves crossing the visited set. 
\end{abstract}

\newpage


\section{Introduction}

Let $(X,f)$ be some dynamical system preserving a probability measure $\mu$
and pick an arbitrary Borel set $A\subset X$ of positive $\mu$-measure. 
By Poincar{\'e}'s recurrence Theorem, $\mu$-almost every point $x\in A$ 
comes back infinitely many times in $A$ upon iterations of $f$.
We denote by $\r_A^n(x)$, $n\geq 1$, the time of the $n$th return of $x\in A$ to $A$.
These times are defined by induction in the following way :
\begin{eqnarray*}
\r_A^1(x)&\pardef&\min\{k\geq 1,\ f^k(x)\in A\},\\
\r_A^{n+1}(x)&\pardef& \r_A^n(x)+\r_A^1(f^{\r_A^n(x)}(x)).
\end{eqnarray*}
(For convenience we set $r_A^0\pardef 0$.)
Assuming that $\mu$ is ergodic, the well-known \kac's formula 
(virtually found in any textbook on ergodic theory) tells that $\int \r_A^1\ d\mu =1$.
We can also introduce the Poincar{\'e} cycles of $A$ with respect to $x\in A$
by setting $\tau_A^{n+1}\pardef \r_A^{n+1}- \r_A^n$, $n\geq 0$.
We obviously have $\r_A^{n}= \sum_{j=1}^{n} \tau_A^j$.
In a note \cite{birkhoff}, Birkhoff showed that
\begin{equation}\label{equ-Kac}
\lim_\ninfty\frac{\r^n_A(x)}{n}=
\frac{1}{\mu(A)}\quad\textup{for}\;\; \mu-\textup{almost every}\;
x\in A\,.
\end{equation}
This note was so overshadowed by his subsequent proof of the Ergodic Theorem
for any integrable function that it escaped notice. The `modern' proof of \eqref{equ-Kac} can
be found for instance in \cite{ka,moy} (where Birkhoff's note is not
cited) and it is a simple consequence of the Ergodic Theorem.

It is natural to study the fluctuations of the convergence in \eqref{equ-Kac}.
We can first ask for large deviations (that is of order $\mathcal{O}(1)$) of this Birkhoff
average of Poincar{\'e} cycles. One can also asks for log-normal fluctuations
(that is of order $\mathcal{O}(1/\sqrt{n})$).
As far as we know, such questions have never been investigated.
It is easy to check that the random variables $\tau_A^n$ are
stationary under the conditional probability measure $\mu_A\pardef
\mu(\cdot\cap A)/\mu(A)$. 
The corresponding process inherits ergodicity from ergodicity of
$(X,f,\mu)$. Set in the language of dynamical systems, this means
that the induced system $(A,f^{\tau_A^1},\mu_A)$ is an ergodic
dynamical system \cite{brown}.
The $\tau_A^n$'s are generally not independent.
Let us mention a notable exception: When $A$ is a state of a countable state Markov chain.
A dynamical realisation of such a process is given by the so-called Gaspard-Wang
map, a piece-wise linear approximation of the Manneville-Pomeau map, see e.g.
\cite{collet}.

The aim of the present article is to prove an accurate large deviation principle
and a central limit theorem for $\r^n_A / n$ when 
$f$ is an Axiom-A diffeomorphism on a Riemannian manifold $M$,
and $\mu$ is the equilibrium state associated to a H{\"o}lder continuous
potential $\varphi$. 
The visited set $A$ will be a Markov rectangle of some basic set (or a finite union
of Markov rectangles).

There are a lot of recent works on return times. In most of them, ``rare
events'' are considered, that is sets $A_n$ such that $\mu(A_n)$ goes to zero as $n$ tends
to infinity. Typically, $A_n$ is a cylinder set and one looks at the rescaled
return times to $A_n$. In many dynamical systems with ``sufficiently strong''
mixing properties, such rescaled returns are shown to be distributed according
to a Poisson law as $n$ goes to infinity (see \cite{AG} for a recent but not
up-to-date review). For a fixed set, the moments of hitting and return
times are studied in \cite{chazottes} in the setting of `strongly mixing' processes.
We would like to emphasize that in the present work the asymptotics
are taken with respect to the {\em number of visits} to some {\em fixed} set $A$.

Our key-result (Theorem \ref{maintheorem}) is the existence of a kind
of ``free energy'' for Poincar{\'e} cycles.
In probabilistic terms, we prove the existence of the scaled-generating cumulant
function associated to $\r^n_A$. We are able to analyse the properties of this
function because we show it is nothing but the logarithm of the largest
eingenvalue of some one-parameter family of transfer operators. These transfer
operators act on the induced system on an unstable leaf of reference crossing
the set $A$ upon consideration. This construction was done in \cite{leplaideur}
for other purposes than studying return times (namely to construct
equilibrium states).
Once we have this free energy for Poincar{\'e} cycles and its properties, we apply
two results (that are little used) to get precise
large deviation estimates and a central limit theorem for Birkhoff averages
of Poincar{\'e} cycles. At a more technical level, let us notice that
since return times are not continuous functions, we cannot apply the
so-called contraction principle of large deviation theory \cite{DZ}
to the empirical measure. We can neither (directly) apply
the known central limit asymptotics which are established for Lipschitz continuous
functions (see \cite{wad} and references therein).

\bigskip

{\it Outline of the article}.
In Section \ref{main}, we state our main result and its consequences. Section
\ref{preparatory} is devoted to some preparatory notions and lemmas.
The proof of the main result is given in Section \ref{proof}. We first
handle the case when $A$ is a single Markov rectangle. We then show
how to extend the result to a finite union of rectangles. 
In Section \ref{corollaries} we derive our large deviation and central
limit theorems.

\section{Statement of results}\label{main}

We refer the reader to the book of Bowen (see \cite{bowen}) for the
precise definitions of Axiom-A diffeomorphisms, equilibrium states and
basic sets.

\bigskip

{\bf Assumptions}.
Throughout we assume that $f$ is a $C^2$ Axiom-A diffeomorphism on a 
compact Riemannian manifold $M$. Let $\Omega$ be a basic set for $f$
and $\varphi:M\to \mathbb{R}$ be a H{\"o}lder continuous function.
We denote by $\mu$ the (unique) equilibrium state
associated to $\varphi$ on $\Omega$.
Finally, $\CR=\{R_i\}$ denotes a finite Markov
partition of $\Omega$ (into more than one rectangle). Let $A\subsetneq
\Omega$ be some finite union
of atoms of the partition $\CR$.

\bigskip

The main result that we are going to prove is the following :

\begin{theorem}\label{maintheorem}
Under the above assumptions, there exists a real number
$\alpha_0=\alpha_0(A)>0$ such that for every $\alpha<\alpha_0$,
$$
\Psi(\alpha)\pardef\lim_{n\rightarrow+\8}\frac1n\log\BBe_{\mu_A}
\left[e^{\alpha \r_A^n}\right]<+\8\ .
$$
Moreover the map $\alpha\mapsto \Psi(\alpha)$ has the following properties:

\noindent 1.  The map $z\mapsto \Psi(z)$ is analytic in a
complex neighborhood of $]-\8,\alpha_0[$.

\noindent 2. It is strictly convex on $]-\8,\alpha_0[$.

\end{theorem}

We shall apply this theorem and a result due to Plachky and Steinebach \cite{PS} to get
precise estimates on large fluctuations on Birkhoff averages
of Poincar{\'e} cycles. 

\begin{theorem}[Large deviations]\label{LD}
Under the above assumptions, we have the following estimates,
for every $u\in (0,\infty)$
$$
\lim_{n\to\infty} \frac{1}{n} \log\ \mu_A \left\{ \frac{\r_A^n}{n}\geq \frac{1}{\mu(A)} +u \right\}
=\inf_{\alpha<\alpha_0} \left\{-\left(\frac{1}{\mu(A)}+u\right) \alpha +  \Psi(\alpha) \right\}
$$
and for every $0<u<1/\mu(A)$
$$
\lim_{n\to\infty} \frac{1}{n} \log\ \mu_A \left\{ \frac{\r_A^n}{n}\leq \frac{1}{\mu(A)} -u \right\}
=\inf_{\alpha<\alpha_0} \left\{-\left(\frac{1}{\mu(A)}-u\right) \alpha +  \Psi(\alpha) \right\}
$$
where $\alpha_0>0$ is the same as in Theorem \ref{maintheorem}.
\end{theorem}

Of course we can replace the $n$th return time, $\r_A^n$, by the Birkhoff
sum of Poincar{\'e} cycles, $\sum_{j=1}^{n} \tau_A^j$, in the previous theorem.

\bigskip

We also obtain a central limit theorem by using Theorem \ref{maintheorem} and
applying a result due to Bryc \cite{bryc}.
Notice that in general it is impossible to deduce a central limit theorem
from a large deviation principle assuming only that the cumulant generating
function is twice differentiable at the origin. Even real-analyticity is not
enough (see a counterexample in \cite{bryc}).

\bigskip

\begin{theorem}[Central limit theorem]\label{CLT}
Under the above assumptions
\begin{equation}\label{convergence}
\lim_{n\to\infty}
\mu_A\left\{
\frac{\r_A^n - n/\mu(A)}{\sigma_A\sqrt{n}} \leq t 
\right\} 
= \frac{1}{\sqrt{2\pi}} \int_{-\infty}^{t} e^{-\frac{\xi^2}{2}}\ d\xi
\end{equation}
where
\begin{equation}\label{variance}
\sigma_A^2 = \Psi''(0)=
\lim_{n\to\infty} \frac1{n} \int \left(\r_A^n- \frac{n}{\mu(A)}
\right)^2  d\mu_A\in\ ]0,+\infty[\,.
\end{equation}
\end{theorem}

We can replace the $n$th return time, $\r_A^n$, by the Birkhoff
sum of Poincar{\'e} cycles, $\sum_{j=1}^{n} \tau_A^j$, in the previous theorem.

\newpage

\noindent {\bf Remarks}

\begin{enumerate}

\item $\alpha_0$ has an explicit expression: it is the difference of
the topological pressures of the system and the topological pressure
of the system obtained by `removing $A$' from the phase space, see formula
\eqref{alpha0} below. 

\item We emphasize that in Theorem \ref{LD} we have limits, not just
liminf and limsup. We could of course formulate the result in terms
the Legendre transform of $\Psi$. We could also consider more general
intervals of large deviations in the spirit of the G{\"a}rtner-Ellis
Theorem \cite{DZ}.

\item 
It is well-known that the values assumed by $\Psi'$ when $\alpha$
ranges from $-\infty$ to $\alpha_0$ give the values of the possible
deviations $u$ around the mean $1/\mu(A)$. It will be easy to check
that $\Psi'(\alpha)$ runs from $0$ to $+\infty$. The fact that
fluctuations above $1/\mu(A)$ can be arbitrary large is due to the
fact that there are points which are typical for (invariant) measures giving
arbitrary small weight to $A$.
Notice that
$\Psi'(0)=\lim_{n\rightarrow+\8}\frac1{n} \BBe_{\mu_A}[\r_A^n]=\frac1{\mu(A)}$.
Indeed, 
$$
\frac1{n} \BBe_{\mu_A}\left(\sum_{j=1}^{n} \tau_A^j\right) = 
\frac1{n} \sum_{j=1}^{n} \BBe_{\mu_A}(\tau_A^j)=
\BBe_{\mu_A}(\tau_A^1)
=\frac1{\mu(A)}\, .
$$
The last equality is \kac formula and the fact that
$\BBe_{\mu_A}(\tau_A^j)=\BBe_{\mu_A}(\tau_A^1)$ 
for all $j\in\mathbb{N}$ is established in \cite{ka}
or \cite{moy}.

\item There is another way to prove a central limit theorem close to
Theorem \ref{CLT}. Let us sketch it when $A$ is a single rectangle.
There is a well-known duality between the $n$th return to $A$ and the number
of occurrences of $A$ up to time $n$. 
Let $N^n_A(x)\pardef \BBone_A(x)+\cdots + \BBone_A(f^{n-1}x)$. This is
the number of visits of the orbits of $x$ to $A$ up to time $n$.
A central limit theorem can be easily derived for $N_A^n$ by using
\cite{PP}. Indeed, we can apply the central limit theorem
given therein to the characteristic function of the one-cylinder
associated to $A$ in the subshift $(\Sigma,\sigma)$ and pull it back to
$(\Omega,f)$.
This means that for all $t\in\R$
$$
\lim_{n\to\infty}
\mu\left\{
\frac{N_A^n - n\mu(A)}{\overline{\sigma}_A\sqrt{n}}\leq t
\right\}
= \frac{1}{\sqrt{2\pi}} \int_{-\infty}^{t} e^{-\frac{\xi^2}{2}}\ d\xi
$$
where
$$
\overline{\sigma}_A^2=
\lim_{n\to\infty}\frac1{n}\textup{Var}(N_A^n)=
\lim_{n\to\infty}
\frac1{n} \int \left(N_A^n -n\mu(A)\right)^2\ d\mu\,.
$$
We know that $0\leq \overline{\sigma}_A^2<\infty$.
Applying Feller's result \cite{feller} we get \eqref{convergence}
with
\begin{equation}\label{variancebis}
\overline{\sigma}_A^2 = \sigma_A^2\ \mu(A)^{3}\, .
\end{equation}
Thus, proving that $\sigma_A^2>0$ is equivalent to
prove that $\overline{\sigma}_A^2>0$ since $A$
is a (non-trivial) Markov rectangle, that is, $0<\mu(A)<1$.
But by \cite{PP} we know that $\overline{\sigma}_A^2=0$
if, and only if, $\BBone_A-\mu(A)$ is a coboundary.
This cannot happen for a Markov rectangle
because this would imply that for any fixed point
$x\in A$, $1=\BBone_A(x)=\mu(A)<1$.
Notice that by following this line of proof, we do not prove that $\sigma_A^2$
is equal to $\Psi'(0)$.

\end{enumerate}

\section{Preparatory lemmas}\label{preparatory}

In this section, we briefly recall the relevant results from \cite{leplaideur}
which are useful for the proof of Theorem \ref{maintheorem} and derive a few
lemmas.

\subsection{Notations}

Let us recall that $\mu$ is the unique equilibrium state
associated to the potential $\varphi$, that is, we have
\begin{equation}\label{equ-equilibre}
h_\mu(f)+\int\varphi\,d\mu=\sup_{\nu}\left(
  h_\nu(f)+\int\varphi\,d\nu\right)
=P_{top}(\varphi,\Omega)
\end{equation}
where the supremum is taken over the set of $f$-invariant probability measures
on $\Omega$. As usual, $h_\nu(f)$ denotes the entropy of the measure $\nu$ and
$P_{top}(\varphi,\Omega)$ the topological pressure on $\Omega$
associated to the potential $\varphi$.
Let $N\geq 2$ be the number of proper rectangles of the Markov partition $\CR$ and
$\CA$ the $N\!\times\! N$-transition matrix defined as
\begin{eqnarray*}
a_{ij}=1\quad \mbox{ if }\quad f^{-1}(\inte{R_j})\cap \inte{R_i}\not=\emptyset\\
a_{ij}=0\quad \mbox{ otherwise.}
\end{eqnarray*}
Let $\Sigma$ be the set of sequences $\ul{x}=\{x_n\}_{n\in\Z}$ such that
for every $n$, $x_n$ belongs to $\{1,\ldots,N\}$ and $a_{x_n x_{n+1}}=1$.
If $\sigma$ denotes the shift map on $\Sigma$, there exists some canonical map
$\pi$ from $\Sigma$ onto $\Omega$ such that the following diagram commutes:
$$
\begin{array}{rcl}
\Sigma&\stackrel{\sigma}{\longrightarrow}&\Sigma\\
\pi\downarrow& \leftturn&\downarrow\pi\\
\Omega&\stackrel{f}{\longrightarrow}&\Omega\\
\end{array}
$$
As the map $\pi$ is also H{\"o}lder continuous, the map $\wt\varphi$ defined by
$$
\wt\varphi\pardef \varphi\circ\pi
$$
is again H{\"o}lder continuous, and there exists a unique equilibrium state
$\wt\mu$ for the dynamical system $(\Sigma,\sigma)$ associated to the
potential $\wt\varphi$.
Moreover $\wt\mu\circ\pi^{-1}=\mu$.

The topological pressure associated to $\wt\varphi$,
$P_{top}(\wt\varphi,\sigma)$, is equal to $P_{top}(\varphi,\Omega)$. 
The cylinder set $[i_0,\ldots,i_n]\subset \Sigma$,
$i_j\in\{1,...,N\}$, $n\geq 0$,
is the set of points $\ul x$ such that $x_j=i_j$ (for every $0\leq
j\leq n$).

\bigskip
Let $g$ be the first return map in $A$:
\begin{eqnarray*}
g:A&\longrightarrow&A \\
x&\longmapsto&f^{\r^1_A(x)}(x).
\end{eqnarray*}
If $x$ is a point in $\Omega$, we denote, as usually, by 
$W^u(x)$, $W^u_{loc}(x)$, $W^s(x)$, $W^s_{loc}(x)$ the unstable
and stable global and local manifolds. 
Local means that the length is equal to some expansive constant, $\eps_0$.
For every set $R$ of small diameter (smaller than $\eps_0$) we set
$$W^i(x,R)\pardef W^i_{loc}(x)\cap R, \mbox{ for }i=u,s.$$
We will assume that the diameter of $\CR$ is smaller than $\eps_0$.

\subsection{The subsystem $(F,g_F)$}

For the sake of definiteness, we set $A=\pi[1]$.
We denote by $F$ some fixed unstable leaf in $A$; namely we
have
$$
F=W^u(x_0,A)
$$
for some fixed point $x_0$ in $\inte {\textup A}$. The system of local
coordinates gives a projection $\pi_F$ from $A$ onto $F$.
This projection is H{\"o}lder continuous. We denote by $g_F$
the map $\pi_F\circ g$. For $x$ in $\Omega$ and $x'$ in
$W^s(x)$, we set
$$
\omega(x,x')=\sum_{k=0}^{+\8}\varphi\circ
f^k(x)-\varphi\circ f^k(x')\ .
$$
The map $\varphi$ is H{\"o}lder continuous, and so, by
contraction on the stable leaves, the previous series
converges. For $x$ in $F$ we set
$\omega(x)=\omega(g(x),g_F(x))$, and
$$
\Phi(x)=\sum_{k=0}^{\r^1_A(x)-1}\varphi\circ
f^k(x)+\omega(x)\ .
$$
This function is defined on a set of full measure with respect to 
any invariant measure.
A simple computation gives the following lemma
\begin{lemme}\label{lem-omega-borne}
There exists some positive constant $C_\omega$ such that
for every $x$ in $F$,
$$\|\omega\|_\8\leq C_\omega.$$
\end{lemme}
The inverses branches of $g_F$ define the family of
$n$-cylinders: for $x$ in $F$ we set
$$C_n(x)\pardef f^{-\r^n_A(x)}(W^u(g^n(x),A)).$$
The $n$-cylinders are well-defined except for the points
in $F$ which do not return infinitely many times in $R$ and
for the points which belong to the orbit of the boundary
$\partial\CR$ of the partition. These two sets of points
will have null-measure for all the measures we are going to
consider. Hence, every $n$-cylinder is a compact set and
the collection of the $n$-cylinders defines a partition of
$F$ (up to the boundary and points which come back less
than $n$ times), which refines the partition into $(n-1)$-cylinders.
An important property is that $g_F^n(C_n(x))=F$ for every
$n$-cylinder. This is a consequence of the Markov property
of the partition. This allows us to define the set of
preimages of some point $x$ in $F$ by $g_F^n$, denoted by
$Pre_n(x)$. Hence, every $n$-cylinder contains exactly one
element of $Pre_n(x)$, for every $x$ in $F$. We define the
Perron-Frobenius-Ruelle operator for any $x$
$$
\CL_S(T)(x)=\sum_{y\in Pre_1(x)}e^{\Phi(y)-\r^1_A(y)S}T(y)
$$
where $S$ is a real parameter. It is proved in
\cite{leplaideur} that there exists a critical value
$S_c$, with $S_c\leq P_{top}(\varphi,\Omega)$, such that
\begin{equation}\label{equ-L_s-cv}
\CL_S(\BBone_F)(x)<+\8
\quad
\textup{for every}\; S> S_c\;\textup{and for every}\;x\in F\,.
\end{equation}
$S_c$ is of course defined as the smallest real number
with this property.
It is also proved that for every $x$ in $F$ and for
$S=P_{top}(\varphi,\Omega)$, we still have
$\CL_S(\BBone_F)(x)<+\8$. A part of the work in the next
subsection will be to prove that in fact $S_c$ is strictly
smaller than the topological pressure
$P_{top}(\varphi,\Omega)$ and indeed equal to the topological pressure
when one removes transitions from or to $A$. 
Here again, the Markov property of the partition and the hyperbolic
structure lead to the next lemma:

\begin{lemme}\label{lem-L_s_cv}
There exists a positive constant $C_\varphi$ which does
not depend on $S$  such that for all $x,y\in F$, $S>S_c$ and integer $n$
$$\frac{1}{C_\varphi}\CL_S^n(\BBone_F)(x)\leq
\CL_S^n(\BBone_F)(y)\leq C_\varphi\CL_S^n(\BBone_F)(x).$$
\end{lemme}

There exists some quasi-metric $\eta$ on $F$ such that for
every $S>S_c$, the operator $\CL_S$ is a quasi-compact
operator on the Banach space $\mathcal{C}_\eta$ of
Lipschitz-continuous function (for the quasi-metric
$\eta$). We recall that the norm $\|\cdot\|_\eta$ on $\mathcal{C}_\eta$
is defined by
$$
\|\phi\|_\eta=\|\phi\|_\8+\sup_{x\not=x'}\frac{|\phi(x)-\phi(x')|}{\eta(x,x')}\,.
$$

The quasi-metric $\eta$ is chosen such that the
$\vartheta$-H{\"o}lder continuous functions on $F$ (where
$\vartheta$ is the H{\"o}lder coefficient of $\varphi$) are
$\eta$-Lipschitz continuous functions. Using the
Ionescu-Tulcea \& Marinescu theorem (see \cite{ionescu}),
we get that there exists some unique probability
measure $m_S$ such that
\begin{equation}\label{def-schauder-tycho}
\CL_S^*(m_S)=\lambda_S m_S
\end{equation}
where $\lambda_S>0$. 
Moreover, there is a unique function $H_S$ in $\mathcal{C}_\eta$ such that
$$
\CL_S H_S = \lambda_S H_S\quad\mbox{ and } \int H_S\,dm_s=1.
$$
This function satisfies
\begin{equation}\label{equ_h_s-bornee}
\frac1{C_0}\leq H_S\leq C_0
\end{equation}
where $C_0>0$. A consequence of Lemma \ref{lem-L_s_cv} is that $C_0$ {\em does not
depend} on $S$.

The measure $\nu_S$ defined by $$d\nu_S=H_S dm_S$$
is the unique equilibrium state for $(F,g_F)$ associated to
$\Phi(\cdot)-S\r^1_A(\cdot)$.
Quasi-compactness of $\CL_S$ means here
that there exist $p=p(S)$ complex numbers of modulus one,
$1=\lambda(1),\ldots,\lambda(p)$ such that
\begin{equation}\label{L_S-quasicompact}
\CL_S=\sum_{i=1}^p\lambda_S\lambda(i)\Psi_i
+\lambda_S\wt\Psi
\end{equation}
where the $\Psi_i$ are linear projectors defined on
$\mathcal{C}_\eta$ with finite rank, $\wt\Psi$ has a spectral radius
strictly smaller than 1, and all the kernels of these
operators contain the images of the others.

A crucial fact is that for $S=P_{top}(\varphi,\Omega)$, we have
$\lambda_S=1$ and $\nu_S$ is the projection on $F$ of the
measure $\mu_A$.

\subsection{Computation of $S_c$.}

We can remove the set $A$ from the Markov partition to define
a new subshift of finite type in the following way.
For the sake of definiteness, we assume that the first line
and the first column of $\CA$ encode the transitions to and
from $A$.
We denote by $\CA'$ the $(N-1)\times (N-1)$ matrix obtained
by removing from $\CA$ the first line and the first column.
$\Sigma'$ will denote the subset in $\Sigma$ of all
sequences $\ul{x}=(x_n)$ such that $a'_{x_nx_{n+1}}=1$. For
convenience we assume that the matrix $\CA'$ is transitive;
if this is not the case, we can restrict our work to the
classes of recurrence (for each class the map $\sigma$
satisfies expansiveness and specification). The map
$\wt\varphi$ can be restricted to $\Sigma'$, thus the
dynamical system $(\Sigma',\sigma)$ admits exactly one
equilibrium state, $\wt\mu'$, with topological pressure
$P_{top}(\wt\varphi,\Sigma')$.
\begin{lemme}\label{trou-de-pression}
With the previous notations,
$P_{top}(\wt\varphi,\Sigma')<P_{top}(\wt\varphi,\Sigma)$.
\end{lemme}

The proof can be found in \cite{CGL}. 
Another proof based on relative entropy can be found in \cite{CFL}.

We can now prove the main lemma of this subsection:

\begin{lemme}\label{prop-S_c-majo}
The critical value $S_c$, defined in \eqref{equ-L_s-cv},
is equal to $P_{top}(\wt\varphi,\Sigma')$. 
\end{lemme}

\begin{proof}
To ease notations we set $P'=P_{top}(\wt\varphi,\Sigma')$ throughout this proof.
Let $x$ be in $F$. By lemma \ref{lem-L_s_cv} we have just to prove
that for every $S>P'$, $\CL_S(\BBone_F)(x)<+\8$
and for $S=P'$, $\CL_S(\BBone_F)(x)=+\8$.  We can choose $x$
such that it does not belong to the set $\bigcup_\Z f^{-n}\partial
\CR$.
Therefore $\pi_F^{-1}(x)$ is a single point $\ul x$ in $\Sigma$; this
also holds for every $y$ in $Pre_1(x)$. 

Let $n$ be some integer. The set, $Pre_1^n(x)$, of points $y$ in
$Pre_1(x)$ such that $\r^1_A(y)=n$
(which is well defined because $f^k(x)$ never belongs to the boundary
of the Markov partition)
is a $(\eps_0,n)$-separated set of points: all these points belong to
$F\subset W^u_{loc}(x)$
and all their images by $f^n$ belong to $W^s_{loc}(x)$. If we also
assume that $\eps_0$
is small enough, the same argument proves that the set $f(Pre_1^n(x))$
is $(\eps_0,(n-2))$-separated (but not necessarily maximal).
Therefore, it follows from \cite{Ru} that there exists some constant
$\kappa_1$, independent of $S$ and $n$, such that 
\begin{equation}\label{equ-patch1.1}
\sum_{y\in Pre_1^n(x)}e^{\Phi(y)}e^{-nS}\leq e^{\kappa_1+C_{\!\omega}}e^{(n-2)P'-nS}.
\end{equation}
Because the matrix $\CA'$ is transitive, there exists some integer $K$
such that ${\CA'}^K$
has only positive entries. We can assume that $n$ is strictly greater
than $2K+2$.
Let us consider the set $E'_{n-2K-2}$ of all words in $\Sigma'$ with
length $n+1-2K-2$.
We pick in each element of $E'_{n-2K-2}$ exactly one point (of
$\Sigma'$);
this collection of points is denoted by $I_n$. As before, there exists
some constant $\kappa_2$ (independent of $n$ and $S$) such that 
\begin{equation}\label{equ-patch1.2}
\sum_{\ul{y}\in I_n}e^{S_{n-2K-1}(\wt\varphi)(\ul{y})}\geq e^{\kappa_2+(n-2K-1)P'},
\end{equation}
where  $S_{n-2K-1}(\wt\varphi)(\cdot)$ means
$\wt\varphi(\cdot)+\wt\varphi\circ\sigma(\cdot)+\cdots+\wt\varphi\circ\sigma^{n-2K-2}(\cdot)$. 

We are now going to prove that points $\ul{y}$ in $I_n$ can be
chosen in such a way that all the
$f^{-K-1}(\pi(\ul{y}))=\pi(\sigma^{-K-1}(\ul{y}))$ are in
$Pre_1^n(x)$.

We pick some $z$ in $Pre_1^n(x)$ and denote by $\ul{z}$ its preimage
by $\pi$;
we also set $\ul{z}=(z_n)_{n \in \Z}$. Hence $z_1$  and $z_{n-1}$ must
be different from 1. Because of our choice of $K$, for each
$W_{n-2K-1}$ in $E'_{n-2K-2}$,
there exist two words $w_K$ and $w'_K$ of length $K$ in $\Sigma'$ such that:
\begin{itemize}
\item the first letter of $w_K$ is $z_1$,

\item the last letter of $w'_K$ is $z_{n-1}$,

\item the word $(w_K W_{n-2K-1} w'_K)$ of length $n-1$ is admissible for $\Sigma'$.
\end{itemize}
Therefore the point
$$
\ul{Z}(W_{n-2K-1})\pardef\ldots,z_{-1},z_0,(w_K,W_{n-2K-1},w'_K),z_n,z_{n+1},\ldots
$$
(with the initial position in $z_0$) belongs to $\Sigma$,
$\pi(\ul{Z}(W_{n-2K-1}))$ is in $Pre_1^n(x)$
and $\sigma^{K+1}(\ul{Z}(W_{n-2K-1}))$ is in the cylinder
$W_{n-2K-1}$.
The sum
$$
\sum_{y\in Pre_1^n(x)}e^{\Phi(y)}e^{-nS}
$$
is greater
than the same sum but restricted to the $\pi(\ul{Z}(W_{n-2K-1}))$'s,
thus by (\ref{equ-patch1.2}) proves that there exists some constant
$\kappa_3$, which does not depend on $S$ and $n$, such that 
\begin{equation}\label{equ-patch1.3}
\sum_{y\in Pre_1^n(x)}e^{\Phi(y)}e^{-nS}\geq e^{\kappa_3-C_\omega}e^{(n-2K-1)P'-nS}.
\end{equation}
As 
$\disp\CL_S(\BBone_F)(x)=\sum_n\sum_{y\in
  Pre_1^n(x)}e^{\Phi(y)}e^{-nS}$
, (\ref{equ-patch1.1}) and (\ref{equ-patch1.3})
prove that $\CL_S(\BBone_F)(x)$ converges for every $S>P'$ and diverges to $+\8$ for every $S\leq P'$.
The lemma is thus proved.
\end{proof}

\section{Proof of Theorem \protect\ref{maintheorem}}\label{proof}

Let us start with the case when $A$ is a single rectangle
of the Markov partition. 
At the end of this section, we briefly
sketch the modifications that are necessary to handle the case when $A$ is a finite
union of rectangles.

Let us set
\begin{equation}\label{alpha0}
\alpha_0=P_{top}(\wt\varphi,\Sigma)-S_c=P_{top}(\wt\varphi,\Sigma)-P_{top}(\wt\varphi,\Sigma')\,.
\end{equation}
Let $\alpha$ be in
$]0,\alpha_0[$. For convenience we will write $P$ in
subscript instead of $P_{top}(\varphi,\Omega)$.
We start by observing that if we pick any point on the
reference unstable leaf $F$, then all points lying in its
stable leaf have the same return times as this point.
That is why we can reduce the problem to the unstable leaf $F$.
We have
\begin{eqnarray}
\nonumber
\BBe_{\mu_A}\left[e^{\alpha \r^n_A}\right]&=&\int_F e^{\alpha
r^n(x)}H_P(x)dm_P(x)
=\int_F\CL_P^n(e^{\alpha \r^n_A}H_P)(x)dm_P(x)\\
&=&\int_F\sum_{y\in
Pre_n(x)}\!\!\!e^{S_{\r^n_A(y)}(\varphi)(y)+\omega(g^n(y))-(P-\alpha)\r^n_A(y)}H_P(y)\,
dm_P(x)\,.\nonumber \\
& & \label{equ1-lim-pgd}
\end{eqnarray}
We first remark that for $\alpha<\alpha_0$ we have
$P-\alpha>S_c$. Hence, using lemmas \ref{lem-omega-borne}-\ref{lem-L_s_cv},
the Markov property and the properties of the function $H_P$ one obtains the existence of
some constant $C$, which depends only on $\varphi$ and $f$
such that for every $x$ in $F$,
\begin{equation}\label{equ2-lim-pgd}
e^{-C}\CL_{P-\alpha}^n(\BBone_F)(x)\leq
\end{equation}
$$
\int_F\sum_{y\in
Pre_n(x)}e^{S_{\r^n_A(y)}(\varphi)(y)+\omega(g^n(y))-(P-\alpha)\r^n_A(y)}H_P(y)\,
dm_P(x) \leq
$$
$$
e^{C}\CL_{P-\alpha}^n(\BBone_F)(x).
$$
If we integrate this double inequality
with respect to the measure
$m_{P-\alpha}$, we obtain the following estimate:
\begin{equation}\label{equ3-lim-pgd}
\log\lambda_{P-\alpha}-\frac{C}{n}\leq
\frac1n\log\BBe_{\mu_A}\left[e^{\alpha
\r^n_A}\right]\leq \log\lambda_{P-\alpha}+\frac{C}{n}\, .
\end{equation}
This proves that for every $\alpha<\alpha_0$
$$
\lim_\ninfty
\frac1n\log\BBe_{\mu_A}\left[e^{\alpha\r^n_A}\right]=
\log\lambda_{P-\alpha}<+\8\,.
$$

\bigskip
\noindent
Let us set $\Psi(\alpha)\pardef\log\lambda_{P-\alpha}$. We now
have to prove that the function $\Psi$ is analytic in some
complex neighborhood of $]-\infty,\alpha_0[$. 
Analyticity of $\Psi$ is equivalent to the
analyticity of the map $S\mapsto \log\lambda_S$ 
in some suitable neighborhood. 
For that purpose, we use a theorem of perturbations
due to Hennion and Herv{\'e} \cite{Hennion-Herve}.

We first have to check that $\CL_S$ has $p=p(S)$
dominating simple eigenvalues (see \cite{Hennion-Herve},
III.2). This is a consequence of ergodicity of the system
$(F,g_F,\nu_S)$. In fact, the proof of the proposition 4.11
in \cite{BDP} can be adapted to our case; hence the
projectors $\Psi_i$ in (\ref{L_S-quasicompact}) have rank
one, and the $\lambda(i)$'s all satisfy
$$\lambda(i)^p=1\ .$$
Thus the operator $\CL_S$ has $p$ dominating simple
eigenvalues.

\begin{lemme}\label{lem-ope-anal}
Let $\CO\CP_{\mathcal{C}_\eta}$ denote the set of linear bounded
operators on $\mathcal{C}_\eta$. Then the map $z\mapsto \CL_z$ is
analytic map from $\{z\in\C,\Re(z)>S_c\}$ to $\CO\CP_{\mathcal{C}_\eta}$.
\end{lemme}

\begin{proof}
For any $z$ in $\C$ with $\Re(z)>S_c$, we set (by extension):
$$
\CL_z(\phi)(x)=\sum_{y\in
Pre_1(x)}e^{\Phi(y)-\r^1_A(y)z}\phi(y)
$$
which  can also be written
$$
\sum_{n=1}^{+\8}\left(\sum_{y\in
Pre_1(x),\ \r^1_A(y)=n}e^{\Phi(y)}\phi(y)\right)e^{-nz}\, .
$$
Let us introduce
$$
K_m(z)(\phi)(x)= \sum_{n=1}^{m}\left(\sum_{y\in
Pre_1(x),\ \r^1_A(y)=n}e^{\Phi(y)}\phi(y)\right)e^{-nz} .
$$
We are going to prove that the sequence $(K_m)(\cdot)$ converges
to $\CL_{.}$, when $m$ goes to $+\8$, as analytic functions from $\Re(z)>S_c$ to
$\CO\CP_{\mathcal{C}_\eta}$.
Let us fix some compact set $\Gamma$ in $\Re(z)>S_c$ and
pick $z$ in $\Gamma$. Let $\phi$ be some function in
$\mathcal{C}_\eta$, with $\|\phi\|_\eta=1$. We want to compute
$$\|K_m(z)(\phi)-\CL_z(\phi)\|_\eta.$$
There exists $S>S_c$ such that  for every $z$ in $\Gamma$,
$\Re(z)>S$, which proves that  the series $K_m(z)(\phi)(x)$
is normally convergent (in $z$) and uniformly convergent in
$x$ to $\CL_z(\phi)(x)$ when $m\to\8$. Hence $K_m(z)$ is uniformly
convergent to $\CL_z$ in $\Gamma$ for the norm $\|\cdot\|_\8$
on $\mathcal{C}_\eta$. 
Let $y$ and $y'$ be two points in $F$ with $y'\in C_1(y)$.
Then we have
$$e^{\Phi(y)}\phi(y)-e^{\Phi(y')}\phi(y')=e^{\Phi(y)}\phi(y)-e^{\Phi(y)}\phi(y')
+e^{\Phi(y)}\phi(y')-e^{\Phi(y')}\phi(y').$$

Therefore, the Lipschitz  properties of the functions
$\phi$ and $\Phi$ (for the quasi-metric $\eta$) and the
expansion on the unstable leaves imply that there exists
some constant $C$ which depends only on $f$ and $\phi$ such
that for every $n$, and for every $x$ and $x'$ in $F$,
\begin{equation}\label{equ-rab1-preuveth}
\left|\sum_{y\in Pre_1(x),\ \r^1_A(y)=n}\left(
e^{\Phi(y)}\phi(y)- e^{\Phi(y')}\phi(y')\right)\;
\right|\leq C\eta(x,x'),
\end{equation}
where for each $y$ in $Pre_1(x)$, $y'$ is the preimage of
$x'$ in the 1-cylinder $C_1(y)$. Inequality
(\ref{equ-rab1-preuveth}) and convergence for the norm
$\|\cdot\|_\8$ imply that the series $K_m(z)$ is uniformly
convergent to $\CL_z$ in $\Gamma$ for the norm
$\|\cdot\|_\eta$. This proves that $z\mapsto \CL_z$ is
analytic.
\end{proof}
By Theorem III.8 from \cite{Hennion-Herve}, for every
$S>S_c$, there exists some open disc in $\C$ centered at
$S$, such that for every $z$ in this disk, $\CL_z$
has $p(S)$ dominating simple eigenvalues,
$\lambda^1(z),\ldots,\lambda^p(z)$ and the maps $z\mapsto
\lambda^i(z)$ are analytic.
Because the map $\log$ is analytic in $\C\setminus\R_-$, we
can conclude that the map $z\mapsto \log\lambda_z$ is analytic
in a complex neighborhood of $]S_c, +\infty[$.

\bigskip

It remains to prove strict convexity of $\alpha\mapsto \Psi(\alpha)$ on 
$]-\infty,\alpha_0[$. Convexity is obvious (by H{\"o}lder inequality).
Strict convexity is equivalent to strict convexity of the map $S\mapsto\log\lambda_S$ 
on $]S_c,+\infty[$.
We will assume that the map $S\mapsto \log\lambda_S$ is not strictly
convex and will arrive to a contradiction
with the fact that $A$ is a proper rectangle of the Markov partition.

If the map $S\mapsto \lambda_S$ is not strictly convex, this means
that its graph contains a straight line interval, say
$\mathcal{I}\subset\, ]S_c,\infty[$.
In turn, this means that
$\lambda'_S/\lambda_S$ 
is constant on $\mathcal{I}$, which means
that $\lambda_S = \Lambda_A e^{\gamma_A S}$ where $\Lambda_{\!A}$, $\gamma_A$ are real constants
{\it a priori} depending on $A$.
Now we invoke the unicity of analytic continuation of analytic functions
to deduce that
$$
\forall z\in\mathbb{C}\quad\textup{with}\quad\Re(z)>S_c,\quad \lambda_z = \Lambda_A e^{-\gamma_A z}\,.
$$
Let us define the smallest return time in $A$ as follows:
\begin{equation}\label{tauA}
\tau(A)=\inf\{k\geq 1: f^{-k}A\cap A\neq \emptyset\}=\inf\{\r^1_A(x):x\in A\}\,.
\end{equation}
Using Lemma \ref{lem-L_s_cv}, we readily get for every $S>S_c$
$$
C_\varphi^{-1} \sum_{n\geq \tau(A)}
\left(\sum_{y\in Pre_1(x), \r_A^1(x)=n} e^{\Phi(y)}\BBone_F(y)\right)
e^{-nS}
\leq
$$
$$
\Lambda_A e^{\gamma_A S}
\leq
$$
$$
C_\varphi \sum_{n\geq \tau(A)}
\left(\sum_{y\in Pre_1(x), \r_A^1(x)=n} e^{\Phi(y)}\BBone_F(y)\right)
e^{-nS}\,.
$$
Letting $S\to\8$ we deduce that $\gamma_A=-\tau(A)$.

Now, observe that
$$
\tau(A)=\Psi'(0)=1/\mu(A)
$$ 
(remember Remark 3 in Section \ref{main});
we use \kac formula to obtain
$$
\tau(A)=\sum_{n=\tau(A)}^{\infty} n \mu_A\{\r_A^1=n \}
\,.
$$
Therefore, $\r^1_A(x)=\tau(A)$ for $\mu_A$-almost every $x$.
Since $\mu_A({\inte {\textup A}})=1$, the topological mixing property imposes that
$\tau(A)$ must equal to one, which is absurd since the Markov partition is
made of more than two rectangles.

Therefore, the function $S\mapsto \log\lambda_S$ is strictly convex on 
$]S_c,\infty[$, so is the function $\alpha\mapsto\Psi(\alpha)=\log\lambda_{P-\alpha}$
on $]-\infty,P-S_c[\ =\ ]-\infty,\alpha_0[\ $.

\bigskip

\noindent {\bf Extension to a finite union of rectangles}

\bigskip

We sketch how to extend the previous proof when
the set $A$ is a union $R_{i_1}\cup\ldots R_{i_k}$ of
$k\geq 2$ rectangles (of the partition $\CR$). 
In each rectangle $R_{i_j}$ we pick some unstable leaf $F_j$, and define
$F$ as the union of the $F_j$'s. Up to the boundary of $\partial\CR$, this
union is a disjoint union. Thus, if we denote by $\pi_j$ the projection
from $R_{i_j}$ onto $F_j$ (namely $[.,F_j]$, see \cite{bowen}).
This defines a map $\pi_F$ from $A$ onto $F\pardef\bigcup F_j$ 
(up to some set with zero-measure for all the measures we are going to consider).
Hence the map $g_F\pardef\pi_F\circ g$ is well defined; the Markov property
of the partition $\CR$ allows again us to define the partitions of
$n$-cylinders; but in this case we will not have $g^n_F(C_n(x))=F$, but only
$g^n_F(C_n(x))=F_j$ for some $j$.

The definition of the operator $\CL_S$ is the same, but the value of the
critical $S$, $S_c$, is changing. 
Let $\CA''$ be the matrix obtained from $\CA$ by removing the lines and the
columns corresponding to the subscripts $i_j$'s. Let $\Sigma''$ be the set
of sequences $\ul{x}$ such that for every $n$, $a''_{x_n,x_{n+1}}=1$ and
let $P''$ be the topological pressure for the dynamical system
$(\Sigma'',\sigma)$ associated to the potential $\wt\varphi$. Then,
$S_c\leq P''<P_{top}(\varphi,\Omega)$.

Lemma \ref{lem-L_s_cv} still holds but only if we pick $x$ and $y$ in the same
$F_j$. Accordingly, the proof of proposition \ref{prop-S_c-majo} can be adapted.

Now, observe that formulas \eqref{equ1-lim-pgd} and (\ref{equ2-lim-pgd}) are still valid,
except that we have to split the integrals over the $F_j$'s. Therefore
(\ref{equ3-lim-pgd}) still holds, for some positive constant $C$.

However, it is important to notice that
transitivity on $\Omega$ implies that for every $j$ and $j'$, the set of
points in $F_j$ which returns infinitely many times in $F_{j'}$ is dense.
Moreover, if $x$ is a point in $F_{j'}$ for every $j$, the set of preimages
of $x$ in $F_j$ (for the map $g_F$) is dense.
Therefore exactness of $m_S$ still holds and $\nu_S$ is hence ergodic.
We can thus again apply Theorem III.8 from \cite{Hennion-Herve}
to get analyticity of the map $S\mapsto \log\lambda_S$.

The proof of strict convexity is analogous to the previous case.

\section{Proofs of Theorem \protect\ref{LD} and \protect\ref{CLT}}\label{corollaries}

\noindent {\bf Proof of Theorem \ref{LD}}.
Plachky--Steinebach's result applies once we remind that since $\Psi$ 
is strictly convex and real-analytic on $]-\infty,\alpha_0[$, the
function $\alpha\mapsto\Psi'(\alpha)$ is strictly increasing on that
interval.

We also observe that $\Psi'(\alpha)\to +\infty$ when
$\alpha\to\alpha_0$. This is equivalent to show that
$\lambda'_S/\lambda_S\to\infty$ when $S\to S_c$.
On the other hand, $\Psi'(\alpha)\to 0$
when $\alpha\to -\infty$. This can be showed by using Lemmas
\ref{lem-L_s_cv} and \ref{prop-S_c-majo}, and convexity.

\bigskip

\noindent {\bf Proof of Theorem \ref{CLT}}.
Bryc's Theorem \cite{bryc} applies. In particular, it says that the
variance is equal to $\Psi''(0)$. 
It is easy to get formula \eqref{variance} by differentiating twice
$\Psi$ and using a classical factorisation.
The only point to be proved is that $\sigma_A>0$. 

It is well-known that the variance can be written as follows
$$
\sigma_A^2=
\BBe_A((\tau^1_A)^2)-\frac{1}{\mu(A)^2} +
2 \sum_{j=2}^{\infty}
\left(\BBe_A(\tau^1_A\ \tau_A^j)-\BBe_A(\tau^1_A)\BBe_A(\tau_A^j)\right)\,.
$$

Now we follow the proof of Proposition 4.12, p. 63, in \cite{PP}
(note that Herglotz's Theorem also applies in our setting).
To this end, we just have to check that 
$\BBe_A(\tau^1_A\ \tau_A^j)-\BBe_A(\tau^1_A)\BBe_A(\tau_A^j)$
decreases exponentially fast to $0$. This fact follows from
the $\psi$-mixing property of the induced system (see \cite{CG}):
$$
\left| \BBe_A(\tau^1_A\ \tau_A^j)-\BBe_A(\tau^1_A)\BBe_A(\tau_A^j)\right| \leq 
$$
$$
\sum_{p,q\in \N} pq \ \left| \mu_A(\tau_A^1 =p, \tau_A^j=q) -
\mu_A(\tau_A^1 =p)\ \mu_A(\tau_A^j=q) \right|\leq 
$$
$$
C \left(\sum_{p\in \N} p \mu_A(B_p)\right)^2\ \theta^{j}= C\mu(A)^{-2} \theta^j
$$
where $C>0$, $0<\theta<1$ and $B_p\pardef \{\tau^1_A=p\}$ (we used
\kac formula).

We conclude that $\sigma^2_A=0$ if and only if $\tau^1_A -1/\mu(A)$
is a $L^2(\mu_A)$ coboundary with respect to $g$, the induced map on $A$.
But if $A$ is a single Markov rectangle, this is impossible.
Indeed, there exists a fixed point for $g$ which is periodic
with period $\tau(A)$ for $f$ ($\tau(A)$ is defined in \eqref{tauA}). 
Reasoning as above, this leads to a contradiction with the fact that
$A$ is a strict subset of $\Omega$ in measure. The case when $A$ is a (finite) union of rectangles
is left to the reader.

Therefore we arrive at the conclusion
that the variance $\sigma_A^2$ defined in \eqref{variance}
is strictly positive. Theorem \ref{CLT} has now a complete
proof.


\end{document}